\documentclass[12pt]{article}
\usepackage{amsmath}
\usepackage{amssymb}
\usepackage{amsthm}
\usepackage{latexsym}
\usepackage{graphicx}

\begin{document}
\theoremstyle{plain}
\newtheorem{theorem}{Theorem}
\newtheorem{prop}[theorem]{Proposition}
\newtheorem{cor}[theorem]{Corollary}
\newtheorem{lemma}[theorem]{Lemma}
\newtheorem{question}[theorem]{Question}
\newtheorem{conj}[theorem]{Conjecture}
\newtheorem{assumption}[theorem]{Assumption}

\theoremstyle{definition}
\newtheorem{definition}[theorem]{Definition}
\newtheorem{notation}[theorem]{Notation}
\newtheorem{condition}[theorem]{Condition}
\newtheorem{example}[theorem]{Example}
\newtheorem{introduction}[theorem]{Introduction}

\numberwithin{theorem}{section}

\makeatletter                          
\let\c@equation\c@theorem              
\makeatother                           
\renewcommand{\theequation}{\arabic{section}.\arabic{equation}}

\makeatletter                          
\let\c@figure\c@theorem              
\makeatother                           
\renewcommand{\thefigure}{\arabic{section}.\arabic{figure}}

\newcommand{\todo}[1]{\vspace{5 mm}\par \noindent
   \marginpar{\textsc{ToDo}}\framebox{\begin{minipage}[c]{0.95 \textwidth}
   \tt #1 \end{minipage}}\vspace{5 mm}\par}

\providecommand{\abs}[1]{\lvert#1\rvert}
\providecommand{\norm}[1]{\lVert#1\rVert}
\providecommand{\Z}{\mathbb{Z}} \providecommand{\R}{\mathbb{R}}
\providecommand{\N}{\mathbb{N}} \providecommand{\C}{\mathbb{C}}
\providecommand{\Q}{{\mathbb{Q}}} \providecommand{\x}{\mathbf{x}}
\providecommand{\y}{\mathbf{y}} \providecommand{\z}{\mathbf{z}}
\providecommand{\boxend}{\hspace{\stretch{1}}$\Box$\\ \ \\}
\long\def\symbolfootnote[#1]#2{\begingroup%
\def\thefootnote{\fnsymbol{footnote}}\footnote[#1]{#2}\endgroup}
\renewcommand{\mod}[1]{\,(\text{mod }#1)}

\providecommand{\D}{\mathcal{D}}

\title{Computing the period of an Ehrhart quasi-polynomial}
\author{Kevin M. Woods\footnote{Partially supported by a Clay Liftoff Fellowship and NSF Grant DMS 0402148.}}

\maketitle

\begin{abstract}
If $P\subset \R^d$ is a rational polytope, then $i_P(t):=\#(tP\cap
\Z^d)$ is a quasi-polynomial in $t$, called the Ehrhart
quasi-polynomial of $P$.  A period of $i_P(t)$ is $\D(P)$, the
smallest $\D\in \Z_+$ such that $\D\cdot P$ has integral vertices.
Often, $\D(P)$ is the minimum period of $i_P(t)$, but, in several
interesting examples, the minimum period is smaller. We prove that,
for fixed $d$, there is a polynomial time algorithm which, given a
rational polytope $P\subset\R^d$ and an integer $n$, decides whether
$n$ is a period of $i_P(t)$.  In particular, there is a polynomial
time algorithm to decide whether $i_P(t)$ is a polynomial.  We
conjecture that, for fixed $d$, there is a polynomial time algorithm
to compute the minimum period of $i_P(t)$. The tools we use are
rational generating functions.
\end{abstract}

\section{Introduction}
Given a rational polytope $P\subset \R^d$ (that is, a bounded subset
of $\R^d$ which is defined by a finite collection of integer linear
inequalities), define the function
\[i_P(t) = \#(tP\cap\Z^d),\]
where $tP$ is $P$ dilated by a factor of $t$. Also define $\D=\D(P)$
to be the smallest $\D\in \Z_+$ such that $\D\cdot P$ has integral
vertices.  Ehrhart proved \cite{Ehrhart62} that $i_P(t)$ is a
quasi-polynomial function with a period of $\D$. In other words,
there exist polynomial functions $f_0(t), f_1(t), \ldots ,
f_{D-1}(t)$, called the \emph{constituents} of $i_P(t)$, such that
\[i_P (t) = f_j(n)\text{ for }t \equiv j \mod{\D}.\]

\begin{example}
$P=[0,\frac{1}{2}]\times [0,\frac{1}{2}]\subset \R^2.$
\end{example}
Then
\[i_P(t)=\left\{%
\begin{array}{ll}
    \left(\frac{t+2}{2}\right)^2, & \hbox{for $t$ even} \\
    \left(\frac{t+1}{2}\right)^2, & \hbox{for $t$ odd} \\
\end{array}%
\right. .\] \boxend

We know that $\D$ is \emph{a} period of the quasi-polynomial
$i_P(t)$.  What is the \emph{minimum} period?  Certainly, it must
divide $\D$.  In most cases, in fact, it is exactly $\D$.  In
certain interesting examples, however, the minimum period is
smaller.

\begin{example}
Given partitions $\lambda$ and $\mu$, define the Gelfand-Tsetlin
polytope $P=P_{\lambda\mu}\subset\R^N$, as in \cite{DM03} (following
the classic \cite{GC50}), where $N$ is defined in terms of the
lengths of $\lambda$ and $\mu$.
\end{example}
Then $\#(P\cap\Z^N)$ is the dimension of the weight $\mu$ subspace
of the irreducible representation of $\mathrm{GL}_n\C$ with highest
weight $\lambda$. Though
\[i_P(t)=\#(P_{t\lambda,t\mu}\cap\Z^N)\] is a polynomial, that is, it has period one (see
\cite{KR86}), $\D(P)$ may be made arbitrarily large by suitable
choice of $\lambda$ and $\mu$ (see \cite{DM03}).

\boxend

\begin{example}
More generally, given partitions $\lambda,$ $\mu,$ and $\nu$ such
that $\abs{\lambda}+\abs{\mu}=\abs{\nu}$, define the \emph{hive
polytope} $P=P_{\lambda\mu}^{\nu}\subset\R^N$ as in \cite{Buch00}
(an exposition of ideas from \cite{KT99}), where $N$ is defined in
terms of the length of $\lambda,$ $\mu,$ and $\nu$.
\end{example}
Then $\#(P\cap\Z^N)$ is the Littlewood-Richardson coefficient
$c_{\lambda\mu}^\nu$, defined to be the multiplicity of $V_{\nu}$
(the highest weight representation of $\mathrm{GL}_n(\C)$
corresponding to $\nu$) in $V_{\lambda}\otimes V_{\mu}$.  Though
\[i_P(t)=\#(P_{t\lambda,t\mu}^{t\nu}\cap\Z^N)\]
is a polynomial (see \cite{DW02}), $\D(P)$ need not be one. \boxend

\begin{example}
Given any $\D\subset\Z_+$ and any $s$ dividing $\D$, let $P$ be the
pentagon with vertices $(0,0), (0,-\frac{1}{s}), (\D,-\frac{1}{s}),
(\D,0),$ and $(1,\frac{\D-1}{\D})$.
\end{example}
In \cite{MW04}, it is shown that this pentagon has $\D(P)=\D$, but
has minimum period $s$. \boxend

These examples raise several questions: When is the minimum period
of $i_p(t)$ less than $\D(P)$?  When is $i_P(t)$ a polynomial?  How
can we tell what the minimum period of $i_P(t)$ is? These questions
are wide-open, though \cite{MW04} gives a geometric characterization
of the polygons $P\subset\R^2$ such that $i_P(t)$ is a polynomial.
Here, we attack these questions from a computational perspective.
Can we find algorithms to answer these questions ``quickly?''

Let us be more precise.  We define the \emph{input size} of an
algorithm to be the number of bits needed to encode the input into
binary.  In particular, the input size of an integer $a$ is
approximately $1+\log_2\abs{a}$ (the number of digits needed to
write $a$ in binary).  An algorithm is called \emph{polynomial time}
if the number of steps it takes is bounded by a certain polynomial
in the input size.  Proving that an algorithm is polynomial time is
generally regarded as proving that it is ``quick,'' at least
theoretically.  See \cite{Papa94} for general background on
algorithms and computation complexity.

Our algorithms will take as input a polytope $P$.  The input size of
a polytope defined by $n$ linear inequalities $\langle
c_i,x\rangle\le b_i$, where $c_i\in\Z^d, b_i\in\Z$, is approximately
\[nd+\sum_{i,j}\log_2\abs{c_{ij}}+\sum_i\log_2\abs{b_i}.\]

We can now state the main theorem, which we will prove in Section 4.

\begin{theorem}
\label{PeriodCheck} Fix $d$.  There is a polynomial time algorithm
which, given a rational polytope $P\subset \R^d$ and an integer
$n>0$, decides whether $n$ is a period of the quasi-polynomial
$i_P(t)$.

In particular, there is a polynomial time algorithm which decides
whether $i_P(t)$ is a polynomial (that is, whether $n=1$ is a
period).
\end{theorem}

\ \\

It is important that we fix $d$ in this theorem, because problems of
this sort become intractable if $d$ is allowed to vary. For example,
the problem of deciding whether $P$ contains an integer point is
NP-hard if $d$ is not fixed.

Na\"{i}vely applying Theorem \ref{PeriodCheck} yields an algorithm
to find the minimum period of $i_P(t)$ which, unfortunately, is not
polynomial time.  We would have to factor $\D(t)$, which would give
us a set of possible $n$, one of which must be the minimal period.
We will prove the following corollary in Section 4.  By a
\emph{polynomial-time reduction}, of Problem A to Problem B, we mean
that, if there was some oracle which could solve Problem B
instantaneously (more precisely, in the amount of time it takes to
output the answer to Problem B), then we could use that oracle to
get a polynomial time algorithm for Problem A.  In other words,
Problem A is ``as easy as'' Problem B.

\begin{cor}
\label{Reduction} Fix $d$.  There is a polynomial-time reduction of
the problem of finding the minimum period of $i_P(t)$, where $P$ is
a $d$-dimensional polytope, to the problem of factoring a natural
number $\D$.
\end{cor}

\ \\

Unfortunately, the problem of factoring is probably hard.  It is not
known to be polynomial time (read, not too hard) or NP-hard (read,
very hard) and is probably somewhere in between.  Nevertheless, we
make the following conjecture.

\begin{conj}
Fix $d$.  There is a polynomial time algorithm, which, given a
$d$-dimensional polytope $P$, computes the minimum period of
$i_P(t)$.
\end{conj}

\ \\

The tools we will use are \emph{rational generating functions}.
Given a set $S\subset\Z^d$, define the generating function
\[f(S;\x)=\sum_{a=(a_1,\ldots,a_d)\in
S}x_1^{a_1}x_2^{a_2}\cdots x_d^{a_d}=\sum_{a\in S}\x^a.\] Sets that
are very large can sometimes be written compactly as rational
generating functions in the form \begin{equation} \label{GFForm}
f(S;\x)=\sum_{i\in I}\alpha_i
\frac{\x^{p_i}}{(1-\x^{b_{i1}})(1-\x^{b_{i2}})\cdots(1-\x^{b_{ik_i}})},
\end{equation}
 where $\x\in \C^d$, $\alpha_i\in\Q$, $p_i\in\Z^d$, and $b_{ij}\in\Z^d\setminus
0$.

\begin{example}
$S=\{0,1,2,\ldots,n\}$, for some $n$.
\end{example}
Then \begin{align*}f(S;x)&=1+x+x^2+\cdots+x^n\\
&=\frac{1-x^{n+1}}{1-x}.
\end{align*}
\boxend

In Section 2, we present several tools to compute and to manipulate
rational generating functions, most of which were proved in either
\cite{BP99} or \cite{BW03}.

Given a rational polytope $P\subset\R^d$, define the generating
function
\[F_P(t,z)=f_0(t)+f_1(t)z+\cdots+f_{\D-1}(t)z^{\D-1},\]
where the $f_i(t)$ are the constituents of $i_P(t)$.  In Section 3,
we will prove the following proposition, which will be useful in the
proof of Theorem \ref{PeriodCheck}.

\begin{prop}
\label{QPGF} Fix $d$.  There is a polynomial time algorithm which,
given a rational polytope $P$, computes $F_P(t,z)$ as a rational
generating function of the form (\ref{GFForm}).
\end{prop}

\ \\

Finally, in Section 4, we prove Theorem \ref{PeriodCheck} and
Corollary \ref{Reduction}.

\section{Rational generating function tools}
In this section, we present several tools to compute and manipulate
rational generating functions.  Except for Lemma \ref{DecideEqual},
which is proved here, they were proved in either \cite{BP99} or
\cite{BW03}.

First we present a tool for creating rational generating functions.

\begin{theorem}\label{BarvPolyhedron} (Theorem 4.4 of \cite{BP99}) Fix $d$.  Then there exists a
polynomial time algorithm which, for any given rational polyhedron
$P\subset\R^d$, computes $f(P\cap\Z^d;\x)$ in the form
\[f(P\cap\Z^d;\x)=\sum_{i\in I}\epsilon_i
\frac{\x^{p_i}}{(1-\x^{a_{i1}})(1-\x^{a_{i2}})\cdots(1-\x^{a_{id}})},\]
where $\epsilon_i\in\{-1,+1\}$, $p_i,a_{ij}\in\Z^d$, and $a_{ij}\ne
0$ for all $i,j$.  In fact, for each $i$,
$a_{i1},a_{i2},\ldots,a_{id}$ is a basis of $\Z^d$.
\end{theorem}

\ \\

\begin{example}
$P$ is the interval $[0,n]$.
\end{example}
Then $P\cap\Z=\{0,1,2,\ldots,n\}$, and we have already computed
$f(P\cap\Z)=\frac{1-x^{n+1}}{1-x}$. \boxend

Once we have computed some rational generating functions, we also
have several tools to manipulate them.

Let $f(\x)$, with $\x\in\C^d$, be a rational function in the form
(\ref{GFForm}), and let $l_1,l_2,\ldots,l_d\in\Z^n$ be integer
vectors.  These vectors define the \emph{monomial map} $\phi:
\C^n\rightarrow \C^d$ given by
\[\z=(z_1,z_2,\ldots,z_n)\mapsto (\z^{l_1},\z^{l_2},\ldots,\z^{l_d}).\]
If the image of $\phi$ does not lie entirely in the poles of
$f(\x)$, we can define the function $g:\C^n\rightarrow \C$ by
\[g(\z)=f\big(\phi(\z)\big),\]
which is regular at almost every point in $\C^n$. Then $g(\z)$ is
$f(\x)$ specialized at $x_i = \z^{l_i}$.  In particular, if $l_i=0$
for all $i$, then $g(\z)$ is $f(1,1,\ldots,1)$.

\begin{example}
$S$ is a finite set. \end{example} Then $f(S;1,1,\ldots,1)=\abs{S}$.
\boxend

We have the following theorem, which states that, given $f(\x)$ as a
short rational generating function, we can find $g(\z)$ quickly.

\begin{theorem}\label{MonomialSub} (Theorem 2.6 of \cite{BW03})
Let us fix $k$, an upper bound on the $k_i$ in (\ref{GFForm}).  Then
there exists a polynomial time algorithm, which, given $f(\x)$ in
the form (\ref{GFForm}) and a monomial map
$\phi:\C^n\rightarrow\C^d$ such that the image of $\phi$ does not
lie entirely in the poles of $f(\x)$, computes
$g(\z)=f\big(\phi(\z)\big)$ in the form
\[g(\z)=\sum_{i\in
I'}\beta_i\frac{\z^{q_i}}{(1-\z^{b_{i1}})(1-\z^{b_{i2}})\cdots(1-\z^{b_{is}})},\]
where $s\le k$, $\beta_i\in\Q$, $q_i,b_{ij}\in\Z^n$, and $b_{ij}\ne
0$ for all $i,j$.
\end{theorem}

\ \\

Now let $g_1(\x)$ and $g_2(\x)$ be Laurent power series given by
\[g_1(\x)=\sum_{m\in\Z^d}\alpha_m\x^m \text{ and }
g_2(\x)=\sum_{m\in\Z^d}\beta_m\x^m.\] Then the Hadamard product
$g=g_1\star g_2$ is defined to be the power series
\[g(\x)=\sum_{m\in\Z^d}\alpha_m\beta_m\x^m.\]

\begin{example} $S_1,S_2$ are subsets of $Z^d$,
\[g_1(\x)=\sum_{m\in S_1}\x^m \text{, and }g_2(\x)=\sum_{m\in S_2}\x^m.\]
\end{example}
Then
\[(g_1\star g_2)(\x)=\sum_{m\in S_1\cap S_2}\x^m.\]
\boxend

More generally, we may take the Hadamard product with respect to a
proper subset of the variables, by defining
\[g_1(\y,\z)\star_z g_2(\y,\z)\]
as above, except with $\alpha_m$ and $\beta_m$ functions of $\y$. We
have the following theorem (which is a slightly more general version
of Lemma 3.4 of \cite{BW03}, but the proof is the same).

\begin{theorem}
\label{LemmaHad} Fix $k$, $d_1$, and $d_2$. Let $\y\in\C^{d_1},
\z\in\C^{d_2}$, and $\x=(\y,\z)$. Then there exists a polynomial
time algorithm which, given $l\in\Z^{d_1+d_2}$ and functions
\begin{align*}
g_1(\x)&=\sum_{i\in I_1}\alpha_i\frac{\x^{p_i}}{(1-\x^{a_{i1}})\cdots(1-\x^{a_{ik}})} \text{ and}\\
g_2(\x)&=\sum_{i\in
I_2}\beta_i\frac{\x^{q_i}}{(1-\x^{b_{i1}})\cdots(1-\x^{b_{ik}})}
\end{align*} such that $\langle l,a_i\rangle, \langle l,b_i\rangle
\ne 0$, computes $g=g_1\star_z g_2$ (where the Laurent power series
are convergent on a neighborhood of
$(e^{l_1},e^{l_2},\ldots,e^{l_d})$).
\end{theorem}

\ \\

Note that $l$ in the input of the algorithm is important.  For
example, if $f(x)=\frac{1}{1-x}$, then $f$ has two possible Laurent
power series expansions
\[f(x)=1+x+x^2+\cdots \text{ and }
f(x)=-x^{-1}-x^{-2}-x^{-3}-\cdots\] convergent on $\abs{x}<1$ and
$\abs{x}>1$, respectively.  In this paper, however, the power series
we examine will actually be Laurent polynomials (which are
convergent on all of $\C^d$), so we will not have to worry about
$l$.

We present one final generating function tool.

\begin{lemma}
\label{DecideEqual} Fix $d$ and $k$.  There is a polynomial time
algorithm which, given rational generating functions $g_1(\x)$ and
$g_2(\x)$ in the form (\ref{GFForm}) which are known to be Laurent
polynomials, decides whether $g_1\equiv g_2$.
\end{lemma}

\emph{Remark: }The lemma is also true if $g_1(\x)$ and $g_2(x)$ are
Laurent power series with an infinite number of terms, but there are
several complications which will be noted in the proof.

{\bf Proof: }Let $h(\x)=g_1(\x)-g_2(\x)$.  We want to decide whether
$h\equiv 0$. Suppose that
\[h(\x)=\sum_{a\in\Z^d}c_a\x^a,\]
and let
\[\tilde{h}(\x)=h(\x)\star h(\x)=\sum_{a\in\Z^d}c^2_a\x^a.\]
We can compute $\tilde{h}$ in polynomial time, using Theorem
\ref{LemmaHad}.  Then $h\equiv 0$ if and only if $\tilde{h}\equiv
0$. Since we know that $h$ is a polynomial, we must simply check
whether $\tilde{h}(1)=\sum_{a\in\Z^d}c^2_a$ is zero, which we can do
in polynomial time using Theorem \ref{MonomialSub}.  If we did not
know that $h$ is polynomial, we would have to be a little more
careful, and here is a sketch of what to do. We can find bounds $M$
such that if $c_a=0$ for all $a$ with $\norm{a}_{\infty}\le M$, then
$h$ is identically zero, using, for example, ideas from Section 5.1
of \cite{WoodsThesis04}. Then if we take the Hadamard product
\[\bar{h}=\tilde{h}\star\left(\frac{x_1^{-M}-x_1^{M+1}}{1-x_1}
\frac{x_2^{-M}-x_2^{M+1}}{1-x_2}\cdots\frac{x_d^{-M}-x_d^{M+1}}{1-x_d}\right),\]
we now have something which is known to be a Laurent polynomial, and
$h$ is identically zero if and only if $\bar{h}(1)=0$. \boxend

\section{Computing the generating function}
{\bf Proof of Proposition \ref{QPGF}: } Computing, say, $f_0(t)$
alone would be easy, by interpolation.  Indeed, first define
\[g_0(s)=f_0(s\D).\]
We may find $g_0(0), g_0(1),\ldots, g_0(d)$ in polynomial time,
using Theorem \ref{BarvPolyhedron}, and then interpolate, as
follows. Let $V$ be the $(d+1)\times (d+1)$ Vandermonde matrix whose
$i,j$ entry is $(i-1)^{j-1}$ as $1\le i,j \le d+1$.  Then, if
$g_0(s)=a_0+a_1s+a_2s^2+\cdots+a_ds^d$, we have the following
equation:
\[V\cdot\begin{bmatrix}a_0 \\ \vdots \\ a_d\end{bmatrix}=\begin{bmatrix}g_0(0)\\ \vdots \\ g_0(d)\end{bmatrix}.\]
Multiplying by the inverse of $V$, we get the coefficients of
$g_0(s)$, and can then easily recover the coefficients of $f_0(t)$.

We cannot, however, do this for each $f_i(t)$, sequentially, in
polynomial time: there are $\D$ of them, and $\D$ may be exponential
in the input size. Instead, we perform all $\D$ interpolations
simultaneously, using generating functions.

For $0\le i\le \D-1$, let
\[g_i(s)=f\left(s\D+i\right).\]
For $0\le j\le d$, let
\[h_j(z)=g_0(j)+g_1(j)z+g_2(j)z^2+\cdots+g_{\D-1}(j)z^{\D-1}.\]
For $0\le i\le \D-1$ and $0\le k\le d$, let $a_{ik}$ be such that
\[g_i(s)=a_{i0}+a_{i1}s+a_{i2}s^2+\cdots+a_{id}s^d,\]
 and let
\[a_k(z)=a_{0k}+a_{1k}z+a_{2k}z^2+\cdots+a_{\D-1,k}z^{\D-1}.\]
Then we have that
\[V\cdot\begin{bmatrix}a_0(z)\\ \vdots \\ a_d(z)\end{bmatrix}=\begin{bmatrix}h_0(z)\\ \vdots \\ h_{d}(z)\end{bmatrix}.\]
Therefore, if we can compute each $h_j(z)$ in polynomial time as
short rational generating functions, then we could compute the
$a_k(z)$ as short rational generating functions by multiplying by
the inverse of $V$.

We compute
\begin{align*}
h_j(z)&=g_0(j)+g_1(j)z+g_2(j)z^2+\cdots+g_{\D-1}(j)z^{\D-1}\\
&=f_0(j\D)+f_1(j\D+1)z+f_2(j\D+2)+\cdots+f_{\D-1}(j\D+\D-1)z^{\D-1}\\
&=i_P(j\D)+i_P(j\D+1)z+i_P(j\D+2)z^2+\cdots+i_P(j\D+\D-1)z^{\D-1},
\end{align*}
as follows.  Given $j$, define the polyhedron
\[Q_j=\big\{(z,\y):\ 0\le z\le \D-1 \text{ and }y\in(j\D+z)P\big\}.\]
Then
\[f(Q_j;z,\y)=\sum_{0\le a\le \D-1}z^a\sum_{b\in (j\D+a)P}y^b,\]
and
\[h_j(z)=f(Q_j;z,1).\]
We may compute $f(Q_j;z,\y)$ in polynomial time, using Theorem
\ref{BarvPolyhedron}, and then perform the substitution $\y=1$,
using Theorem \ref{MonomialSub}.

We have shown that we can construct the generating functions
$a_k(z)$, for $1\le k\le d$, in polynomial time.  We must now use
these generating functions to compute
\[F_P(t,z)=f_0(t)+f_1(t)z+\cdots+f_{\D-1}(t)z^{\D-1}.\]

Since, for $0\le j \le \D-1$,
\[g_j(s)=a_{j0}+a_{j1}s+\cdots+a_{jd}s^d\]
and
\[f_j(t)=g_j\left(\frac{t-j}{\D}\right),\]
we have that
\[f_j(t)=a_{j0}+a_{j1}\frac{t-j}{\D}+\cdots+a_{jd}\left(\frac{t-j}{\D}\right)^d\]
and
\[F_P(t,z)=\left\{\begin{array}{cccccccc}
 &a_{00}&+&a_{01}\frac{t}{\D}&+&\cdots&+&a_{0d}\left(\frac{t}{\D}\right)^d\\
+&a_{10}z&+&a_{11}\frac{t-1}{\D}z&+&\cdots&+&a_{1d}\left(\frac{t-1}{\D}\right)^dz\\
 &\vdots & & \vdots& & & & \vdots\\
+&a_{\D-1,0}z^{\D-1}&+&a_{\D-1,1}\frac{t-\D+1}{\D}z^{\D-1}&+&\cdots&+&a_{\D-1,d}\left(\frac{t-\D+1}{\D}\right)^dz^{\D-1}
\end{array}\right. .
\]

For $0\le k \le d$, define
\[b_k(t,z)=a_{0k}\left(\frac{t}{\D}\right)^k+a_{1k}\left(\frac{t-1}{\D}\right)^kz+\cdots+a_{\D-1,k}\left(\frac{t-\D+1}{\D}\right)^kz^{\D-1}.\]
Then
\[F_P(t,z)=b_0(t,z)+b_1(t,z)+\cdots+b_d(t,z).\]
For each $k$, we will compute $b_k(t,z)$ from
\[a_k(z)=a_{0k}+a_{1k}z+a_{2k}z^2+\cdots+a_{\D-1,k}z^{\D-1}.\]
In fact
\[b_k(t,z)=a_k(z)\star_z\left[\left(\frac{t}{\D}\right)^k+\left(\frac{t-1}{\D}\right)^kz+\cdots+\left(\frac{t-\D+1}{\D}\right)^kz^{\D-1}\right],\]
and
\[\left(\frac{t}{\D}\right)^k+\left(\frac{t-1}{\D}\right)^kz+\cdots+\left(\frac{t-\D+1}{\D}\right)^kz^{\D-1}\]
can be computed as a short rational generating function in
polynomial time, by expanding all of the terms and repeatedly using
the fact that, for any $k$, $\sum_{i=1}^{\infty}i^kz^i$ is
$\left(z\frac{d}{dz}\right)^k\left(\frac{1}{1-z}\right)$. Therefore
we can compute the $b_k(t,z)$ and hence $F_P(t,z)$ in polynomial
time. \boxend

\section{Deciding whether $n$ is a period}
{\bf Proof of Theorem \ref{PeriodCheck}: } Given $n$ and $P$, we
want to decide whether $n$ is a period of the quasi-polynomial
$i_P(t)$.  Using Proposition \ref{QPGF}, we may compute the
generating function
\[F_P(t,z)=f_0(t)+f_1(t)z+\cdots+f_{\D-1}(t)z^{\D-1}.\]
Define the generating function
\[G_{n,P}(t,z)=f_n(t)+f_{n+1}(t)z+\cdots+f_{\D-1}(t)z^{\D-n-1}+f_0(t)z^{\D-n}+f_1(t)z^{\D-n+1}+\cdots+f_{n-1}(t)z^{\D-1}.\]
Then $n$ is a period of $i_P(t)$ if and only if $F_P(t,z)\equiv
G_{n,p}(t,z)$.  We must show how to compute $G_{n,P}$ in polynomial
time.  Note that
\begin{align*}F_P(t,z)\star_z \left(\frac{z^n-z^{\D}}{1-z}\right)&=F_P(t,z)\star_z \left(z^n+z^{n+1}+\cdots+z^{\D-1}\right)\\
&=f_n(t)z^n+f_{n+1}(t)z^{n+1}+\cdots+f_{\D-1}(t)z^{\D-1}
\end{align*}
and
\begin{align*}F_P(t,z)\star_z \left(\frac{1-z^{n}}{1-z}\right)&=F_P(t,z)\star_z \left(1+z+\cdots+z^{n-1}\right)\\
&=f_0(t)+f_{1}(t)z+\cdots+f_{n-1}(t)z^{n-1}.
\end{align*}
Then
\[G_{n,P}(t,z)=\left[F_P(t,z)\star_z \left(\frac{z^n-z^{\D}}{1-z}\right)\right]z^{-n}+\left[F_P(t,z)\star_z \left(\frac{1-z^{n}}{1-z}\right)\right]z^{\D-n}.\]
This can be computed in polynomial time, using Theorem
\ref{LemmaHad}.

We can decide whether $F_P(t,z)\equiv G_{n,p}(t,z)$ using Lemma
\ref{DecideEqual}, in polynomial time, and the proof follows.
\boxend

{\bf Proof of Corollary \ref{Reduction}: } Compute $\D=\D(P)$ by
taking the least common multiple of the denominators of all of the
coordinates of the vertices of $P$.  Assume that we can find the
prime factorization of $\D$ using an oracle.  Initialize the
following loop with $n_0:=\D$.
\begin{enumerate}
    \item After the $j$th iteration of the loop, $n_j$ is known to
    be a period of $i_P(t)$.
    \item For each prime factor $p$ of $n_j$, decide whether
    $\frac{n_j}{p}$ is a period of $i_P(t)$.
    \begin{itemize}
        \item If none are periods, then $n_j$ is the minimum period
        of $i_P(t)$, and we are done.
        \item if $\frac{n_j}{p}$ is a period of $i_P(t)$ for some
        $p$, then repeat the process with $n_{j+1}=\frac{n_j}{p}$.
    \end{itemize}
\end{enumerate}
This loop must terminate, because eventually we would have $n_j=1$.
\boxend

\section*{Acknowledgements}
Many thanks to Matthias Beck for helpful conversations.  These
results were originally presented at the Mathematisches
Forschungsinstitut Oberwolfach mini-workshop
``Ehrhart-Quasipolynomials: Algebra, Combinatorics, and Geometry.''

\providecommand{\bysame}{\leavevmode\hbox
to3em{\hrulefill}\thinspace}
\providecommand{\MR}{\relax\ifhmode\unskip\space\fi MR }
\providecommand{\MRhref}[2]{%
  \href{http://www.ams.org/mathscinet-getitem?mr=#1}{#2}
} \providecommand{\href}[2]{#2}

\

\noindent\textsc{Department of Mathematics, University of
California,
Berkeley 94720}\\
\emph{Email: }\texttt{kwoods@math.berkeley.edu}


\begin{thebibliography}{DLM03}

\bibitem[BP99]{BP99}
Alexander Barvinok and James Pommersheim, \emph{An algorithmic
theory of
  lattice points in polyhedra}, New Perspectives in Algebraic Combinatorics
  (Berkeley, CA, 1996--97), Math. Sci. Res. Inst. Publ., vol.~38, Cambridge
  Univ. Press, Cambridge, 1999, pp.~91--147. \MR{2000k:52014}

\bibitem[Buc00]{Buch00}
Anders~Skovsted Buch, \emph{The saturation conjecture (after {A}.\
{K}nutson
  and {T}.\ {T}ao)}, Enseign. Math. (2) \textbf{46} (2000), no.~1-2, 43--60,
  With an appendix by William Fulton. \MR{MR1769536 (2001g:05105)}

\bibitem[BW03]{BW03}
Alexander Barvinok and Kevin Woods, \emph{Short rational generating
functions
  for lattice point problems}, J. Amer. Math. Soc. \textbf{16} (2003), no.~4,
  957--979 (electronic). \MR{1 992 831}

\bibitem[DLM03]{DM03}
Jes\'us De~Loera and Tyrrell McAllister, \emph{Vertices of
{G}elfand-{T}setlin
  polytopes}, preprint, arXiv:math.CO/0309329, 2003.

\bibitem[DW02]{DW02}
Harm Derksen and Jerzy Weyman, \emph{On the
{L}ittlewood-{R}ichardson
  polynomials}, J. Algebra \textbf{255} (2002), no.~2, 247--257. \MR{MR1935497
  (2003i:16021)}

\bibitem[Ehr62]{Ehrhart62}
Eug{\`e}ne Ehrhart, \emph{Sur les poly\`edres rationnels
homoth\'etiques \`a
  {$n$}\ dimensions}, C. R. Acad. Sci. Paris \textbf{254} (1962), 616--618.
  \MR{MR0130860 (24 \#A714)}

\bibitem[GC50]{GC50}
Izrail Gelfand and M.~L. Cetlin, \emph{Finite-dimensional
representations of
  the group of unimodular matrices}, Doklady Akad. Nauk SSSR (N.S.) \textbf{71}
  (1950), 825--828. \MR{MR0035774 (12,9j)}

\bibitem[KR86]{KR86}
Anatoli Kirillov and Nikolai Reshetikhin, \emph{The {B}ethe ansatz
and the
  combinatorics of {Y}oung tableaux}, Zap. Nauchn. Sem. Leningrad. Otdel. Mat.
  Inst. Steklov. (LOMI) \textbf{155} (1986), no.~Differentsialnaya Geometriya,
  Gruppy Li i Mekh. VIII, 65--115, 194. \MR{MR869577 (88i:82020)}

\bibitem[KT99]{KT99}
Allen Knutson and Terence Tao, \emph{The honeycomb model of {${\rm
GL}\sb
  n({\bf C})$} tensor products. {I}. {P}roof of the saturation conjecture}, J.
  Amer. Math. Soc. \textbf{12} (1999), no.~4, 1055--1090. \MR{MR1671451
  (2000c:20066)}

\bibitem[MW04]{MW04}
Tyrrell McAllister and Kevin Woods, \emph{The minimum period of the
{E}hrhart
  quasi-polynomial of a rational polytope}, to appear in \emph{Journal of
  Combinatorial Theory, Series A}, 2004.

\bibitem[Pap94]{Papa94}
Christos Papadimitriou, \emph{Computational {C}omplexity},
Addison-Wesley
  Publishing Company, Reading, MA, 1994. \MR{95f:68082}

\bibitem[Woo04]{WoodsThesis04}
Kevin Woods, \emph{Rational generating functions and lattice point
sets}, Ph.D.
  thesis, University of Michigan, 2004.

\end{thebibliography}
\end{document}